\documentclass[12pt,twoside]{article}
\usepackage{latexsym,amsmath,amssymb,hyperref}
\def\eqvsp{}  \newdimen\paravsp  \paravsp=1.3ex

\def\,{\mskip 3mu} \def\>{\mskip 4mu plus 2mu minus 4mu} \def\;{\mskip 5mu plus 5mu} \def\!{\mskip-3mu}
\def\dispmuskip{\thinmuskip= 3mu plus 0mu minus 2mu \medmuskip=  4mu plus 2mu minus 2mu \thickmuskip=5mu plus 5mu minus 2mu}
\def\textmuskip{\thinmuskip= 0mu                    \medmuskip=  1mu plus 1mu minus 1mu \thickmuskip=2mu plus 3mu minus 1mu}
\textmuskip
\def\beq{\eqvsp\dispmuskip\begin{equation}}    \def\eeq{\eqvsp\end{equation}\textmuskip}
\def\beqn{\eqvsp\dispmuskip\begin{displaymath}}\def\eeqn{\eqvsp\end{displaymath}\textmuskip}
\def\bqa{\eqvsp\dispmuskip\begin{eqnarray}}    \def\eqa{\eqvsp\end{eqnarray}\textmuskip}
\def\bqan{\eqvsp\dispmuskip\begin{eqnarray*}}  \def\eqan{\eqvsp\end{eqnarray*}\textmuskip}

\newtheorem{theorem}{Theorem}

\newtheorem{principle}[theorem]{Principle}
\newenvironment{keywords}{\centerline{\bf\small
Keywords}\begin{quote}\small}{\par\end{quote}\vskip 1ex}

\newtheorem{myexample}[theorem]{Example}
\def\fexample#1#2#3{\vspace{-1ex}\begin{myexample}[#2]\label{#1}\rm #3
\hspace*{\fill} $\diamondsuit\quad$ \end{myexample}\vspace{-2ex} }
\def\paradot#1{\vspace{\paravsp plus 0.5\paravsp minus 0.5\paravsp}\noindent{\bf\boldmath{#1.}}}

\def\req#1{(\ref{#1})}

\def\eps{\varepsilon}

\def\fr#1#2{{\textstyle{#1\over#2}}}

\def\SetR{I\!\!R}
\def\SetN{I\!\!N}

\def\SetZ{Z\!\!\!Z}
\def\qmbox#1{{\quad\mbox{#1}\quad}}
\def\e{{\rm e}}                        

\def\P{{\rm P}}                         

\def\v{\boldsymbol}
\def\trp{{\!\top\!}}

\def\G{\Gamma}
\def\a{\alpha}

\def\g{\gamma}
\def\s{\sigma}

\def\b{\beta}
\def\l{\lambda}
\def\r{\rho}

\def\D{{\cal D}}
\def\F{{\cal F}}
\def\M{{\cal M}}
\def\N{{\cal N}}
\def\R{{\cal R}}
\def\X{{\cal X}}
\def\Y{{\cal Y}}

\def\L{\text{\rm Loss}}
\def\h{\hat}

\def\Rank{\text{\rm Rank}}
\def\LR{\text{\rm LR}}
\def\tr{\text{\rm tr}}
\def\Gauss{\text{Gauss}}
\def\KL{\text{\rm KL}}

\begin{document}

\title{\bf\Large\hrule height5pt \vskip 4mm
The Loss Rank Principle for Model Selection
\vskip 4mm \hrule height2pt}
\author{{\bf Marcus Hutter}\\[3mm]
\normalsize RSISE$\,$@$\,$ANU and SML$\,$@$\,$NICTA \\
\normalsize Canberra, ACT, 0200, Australia \\
\normalsize \texttt{marcus@hutter1.net \ \  www.hutter1.net}
}
\date{February 2006}
\maketitle
\vspace*{-5ex}

\begin{abstract}\noindent
We introduce a new principle for model selection in regression and
classification.
Many regression models are controlled by some smoothness or
flexibility or complexity parameter $c$, e.g.\ the number of
neighbors to be averaged over in k nearest neighbor (kNN) regression
or the polynomial degree in regression with polynomials. Let $\h f_D^c$
be the (best) regressor of complexity $c$ on data $D$.
A more flexible regressor can fit more data $D'$ well than a more
rigid one. If something (here small loss) is easy to achieve it's
typically worth less. We define the loss rank of $\h f_D^c$ as the
number of other (fictitious) data $D'$ that are fitted better by
$\h f_{D'}^c$ than $D$ is fitted by $\h f_D^c$. We suggest selecting the
model complexity $c$ that has minimal loss rank (LoRP).
Unlike most penalized maximum likelihood variants (AIC,BIC,MDL),
LoRP only depends on the regression functions and the loss function.
It works without a stochastic noise model, and is directly
applicable to any non-parametric regressor, like kNN.
In this paper we formalize, discuss, and motivate LoRP, study it for
specific regression problems, in particular linear ones, and compare
it to other model selection schemes.
\def\contentsname{\centering\normalsize Contents}
{\parskip=-2.7ex\tableofcontents}
\end{abstract}

\vspace*{-2ex}

\newpage
\section{Introduction}\label{secIntro}

\paradot{Regression}
Consider a regression or classification problem in which we want to
determine the functional relationship $y_i\approx f_{true}(x_i)$
from data $D = \{(x_1,y_1),...,(x_n,y_n)\}\in\D$, i.e.\ we seek a
function $f_D$ such that $f_D(x)$ is close to the unknown
$f_{true}(x)$ for all $x$. One may define regressor $f_D$ directly, e.g.\
`average the $y$ values of the $k$ nearest neighbors (kNN) of $x$
in $D$', or select the $f$ from a class of functions $\F$ that has
smallest (training) error on $D$. If the class $\F$ is not too
large, e.g. the polynomials of fixed reasonable degree $d$, this
often works well.

\paradot{Model selection}
What remains is to select the right model complexity $c$, like $k$
or $d$. This selection cannot be based on the training error, since
the more complex the model (large $d$, small $k$) the better the fit
on $D$ (perfect for $d=n$ and $k=1$). This problem is called
overfitting, for which various remedies have been suggested:

We will not discuss empirical test set methods like
cross-validation, but only training set based methods. See e.g.\
\cite{MacKay:92} for a comparison of cross-validation with Bayesian
model selection. Training set based model
selection methods allow using all data $D$ for regression. The
most popular ones can be regarded as penalized versions of Maximum
Likelihood (ML). In addition to the function class $\F$, one has to
specify a sampling model $\P(D|f)$, e.g.\ that the $y_i$ have
independent Gaussian distribution with mean $f(x_i)$. ML chooses $\h
f_D^c=\arg\max_{f\in\F_c}\P(D|f)$, Penalized ML (PML) then chooses
$\h c=\arg\min_c\{-\log\P(D|\h f_D^c)+$Penalty$(c)\}$, where the
penalty depends on the used approach (MDL \cite{Rissanen:78}, BIC
\cite{Schwarz:78}, AIC \cite{Akaike:73}). In particular, modern MDL
\cite{Gruenwald:04} has sound exact foundations and works very well
in practice. All PML variants rely on a proper sampling model (which
may be difficult to establish), ignore (or at least do not tell how
to incorporate) a potentially given loss function, and are typically
limited to (semi)parametric models.

\paradot{Main idea}
The main goal of the paper is to establish a criterion for selecting
the ``best'' model complexity $c$ %
based on regressors $\h f_D^c$ given as a black box without insight into the
origin or inner structure of $\h f_D^c$, %
that does not depend on things often not given (like a stochastic noise model), %
and that exploits what is given (like the loss function).
The key observation we exploit is that large classes $\F_c$ or more
flexible regressors $\h f_D^c$ can fit more data $D'\in\D$ well than
more rigid ones, e.g.\ many $D'$ can be fit well with high order
polynomials. We define the {\em loss rank} of $\h f_D^c$ as the number of
other (fictitious) data $D'\in\D$ that are fitted better by
$f_{D'}^c$ than $D$ is fitted by $\h f_D^c$, as measured by some loss
function.
The loss rank is large for regressors fitting $D$ not well {\em and} for
too flexible regressors (in both cases the regressor fits many other
$D'$ better). The loss rank has a minimum for not too flexible
regressors which fit $D$ not too bad. We claim that minimizing
the loss rank is a suitable model selection criterion, since it
trades off the quality of fit with the flexibility of the model.
Unlike PML, our new Loss Rank Principle (LoRP) works without a noise
(stochastic sampling) model, and is directly applicable to any
non-parametric regressor, like kNN.

\paradot{Contents}
In Section \ref{secLoRP}, after giving a brief introduction to
regression, we formally state LoRP for model selection.
To make it applicable to real problems, we have to generalize it
to continuous spaces and regularize infinite loss ranks.
In Section \ref{secLM} we derive explicit expressions for the loss
rank for the important class of linear regressors, which includes
kNN, polynomial, linear basis function (LBFR), Kernel, and
projective regression.
In Section \ref{secBayes} we compare linear LoRP to Bayesian model
selection for linear regression with Gaussian noise and prior, and
in Section \ref{secOMS} to PML, in particular MDL, BIC, AIC, and
MacKay's \cite{MacKay:92} and Hastie's et al.\ \cite{Hastie:01}
trace formulas for the effective dimension.
In this paper we just scratch at the surface of LoRP. Section
\ref{secMisc} contains further considerations, to be elaborated on
in the future.

\section{The Loss Rank Principle (LoRP)}\label{secReg}\label{secLoRP}

After giving a brief introduction to regression, classification,
model selection, overfitting, and some reoccurring examples
(polynomial regression Example \ref{exPR} and kNN Example
\ref{exKNN}), we state our novel Loss Rank Principle for model
selection. We first state it for classification (Principle
\ref{prLRD} for discrete values), and then generalize it for
regression (Principle \ref{prLRC} for continuous values), and exemplify
it on two (over-simplistic) artificial Examples \ref{exSD} and
\ref{exSC}. Thereafter we show how to regularize LoRP for realistic
regression problems.

\paradot{Setup}
We assume data $D = (\v x,\v y) := \{(x_1,y_1),...,(x_n,y_n)\} \in
(\X\times\Y)^n=:\D$ has been observed. We think of the $y$ as having an
approximate functional dependence on $x$, i.e.\ $y_i\approx
f_{true}(x_i)$, where $\approx$ means that the $y_i$ are distorted
by noise or otherwise from the unknown ``true'' values
$f_{true}(x_i)$.

\paradot{Regression and classification}
In regression problems $\Y$ is typically (a subset of) the real
numbers $\SetR$ or some more general measurable space like
$\SetR^m$. In classification, $\Y$ is a finite set or at least
discrete. We impose no restrictions on $\X$. Indeed, $\v x$ will
essentially be fixed and plays only a spectator role, so we will
often notationally suppress dependencies on $\v x$.
The goal of regression is to find a function
$f_D\in\F\subset\X\to\Y$ ``close'' to $f_{true}$ based on the past
observations $D$. Or phrased in another way: we are interested in a
regression function $r:\D\to\F$ such that $\h y:=r(x|D)\equiv
r(D)(x)\equiv f_D(x)\approx f_{true}(x)$ for all $x\in\X$.

\paradot{Notation}
We will write $(x,y)$ or $(x_0,y_0)$ for generic data points, use
vector notation $\v x=(x_1,...,x_n)^\trp$ and $\v y=(y_1,...,y_n)^\trp$,
and $D'=(\v x',\v y')$ for generic (fictitious) data of size $n$.

\fexample{exPR}{polynomial regression}{
For $\X=\Y=\SetR$, consider the set $\F_d:=\{f_{\v w}(x)=w_d
x^{d-1}+...w_2 x+w_1 : \v w\in\SetR^d \}$ of polynomials of degree
$d-1$. Fitting the polynomial to data $D$, e.g.\ by least squares
regression, we estimate $\v w$ with $\v{\h w}_D$. The regression
function $\h y=r_d(x|D)=f_{\v{\h w}_D}(x)$ can be written down in
closed form (see Example \ref{exLBFR}).
} 

\fexample{exKNN}{k nearest neighbors, kNN}{
Let $\Y$ be some vector space like $\SetR$ and $\X$ be a metric
space like $\SetR^m$ with some (e.g.\ Euclidian) metric
$d(\cdot,\cdot)$. kNN estimates $f_{true}(x)$ by averaging the
$y$ values of the $k$ nearest neighbors ${\N}_k(x)$ of $x$ in $D$, i.e.\
$r_k(x|D)={1\over k}\sum_{i\in\N_k(x)} y_i$ %
with $|\N_k(x)|=k$ such that $d(x,x_i)\leq d(x,x_j)$ %
for all $i\in\N_k(x)$ and $j\not\in\N_k(x)$.
} 

\paradot{Parametric versus non-parametric regression}
Polynomial regression is an example of parametric regression in the
sense that $r_d(D)$ is the optimal function from a family of
functions $\F_d$ indexed by $d<\infty$ real parameters ($\v w$). In
contrast, the kNN regressor $r_k$ is directly given and is not based
on a finite-dimensional family of functions. In general, $r$ may be
given either directly or be the result of an optimization process.

\paradot{Loss function}
The quality of fit to the data is usually measured by a loss function
$\L(\v y,\v{\h y})$, where $\h y_i=\h f_D(x_i)$ is an estimate of $y_i$.
Often the loss is additive: $\L(\v y,\v{\h y})=\sum_{i=1}^n\L(y_i,\h
y_i)$. If the class $\F$ is not too large, good regressors $r$ can
be found by minimizing the loss w.r.t.\ all $f\in\F$. For instance,
$r_d(D)=\arg\min_{f\in\F_d}\sum_{i=1}^n (y_i-f(x_i))^2$ and $\h
y=r_d(x|D)$ in Example
\ref{exPR}.

\paradot{Regression class and loss}
In the following we assume a (typically countable) class of
regressors $\R$ (whatever their origin), e.g.\ the kNN regressors
$\{r_k:k\in\SetN\}$ or the least squares polynomial
regressors $\{r_d:d\in\SetN_0\}$. Note that unlike $f\in\F$,
regressors $r\in\R$ are not functions of $x$ alone but depend
on all observations $D$, in particular on $\v y$.
Like for functions $f$, we can compute the loss of each regressor
$r\in\R$:
\beqn
  \L_r(D) \;\equiv\; \L_r(\v y|\v x) \;:=\; \L(\v y,\v{\h y})
  \;=\; \sum_{i=1}^n \L(y_i,r(x_i|\v x,\v y))
\eeqn
where $\h y_i=r(x_i|D)$ in the third expression, and the last
expression holds in case of additive loss.

\paradot{Overfitting}
Unfortunately, minimizing $\L_r$ w.r.t.\ $r$ will typically {\em
not} select the ``best'' overall regressor. This is the well-known
overfitting problem. In case of polynomials, the classes
$\F_d\subset\F_{d+1}$ are nested, hence $\L_{r_d}$ is monotone
decreasing in $d$ with $\L_{r_n}\equiv 0$ perfectly fitting the
data. In case of kNN, $\L_{r_k}$ is more or less an increasing
function in $k$ with perfect regression on $D$ for $k=1$, since no
averaging takes place.
In general, $\R$ is often indexed by a ``flexibility'' or smoothness
or complexity parameter, which has to be properly determined. More
flexible $r$ can closer fit the data and hence have smaller loss,
but are not necessarily better, since they have higher variance.
Clearly, too inflexible $r$ also lead to a bad fit (``high bias'').

\paradot{Main goal}
The main goal of the paper is to establish a selection criterion for
the ``best'' regressor $r\in\R$
\begin{itemize}\parskip=0ex\parsep=0ex\itemsep=0ex
\item based on $r$ given as a black box that does not require insight into the
origin or inner structure of $r$,
\item that does not depend on things often not given (like a stochastic noise model),
\item that exploits what is given (like the loss function).
\end{itemize}
While for parametric (e.g.\ polynomial) regression, MDL and Bayesian
methods work well (effectively the number of parameters serve as
complexity penalty), their use is seriously limited for
non-parametric black box $r$ like kNN or if a stochastic/coding
model is hard to establish (see Section \ref{secBayes} for a
detailed comparison).

\paradot{Main idea: loss rank}
The key observation we exploit is that a more flexible $r$ can fit
more data $D'\in\D$ well than a more rigid one. For instance, $r_d$
can perfectly fit all $D'$ for $d=n$, all $D'$ that lie on a
parabola for $d=3$, but only linear $D'$ for $d=2$. We consider
discrete $\Y$ i.e.\ classification first, and fix $\v x$. $\v y$ is
the observed data and $\v y'$ are fictitious others.

Instead of minimizing the unsuitable $\L_r(\v y|\v x)$ w.r.t.\ $r$,
we could ask how many $\v y'\in\Y^n$ lead to smaller $\L_r$ than $\v
y$. Many $\v y'$ have small loss for flexible $r$, and so smallness
of $\L_r$ is less significant than if $\v y$ is among very few other
$\v y'$ with small $\L_r$. We claim that the loss rank of $\v y$
among all $\v y'\in\Y^n$ is a suitable measure of fit. We define the
rank of $\v y$ under $r$ as the number of $\v y'\in\Y^n$ with
smaller or equal loss than $\v y$:
\bqa\label{eqRank}
  \Rank_r(\v y|\v x) \;\equiv\; \Rank_r(L)
  &:=& \#\{\v y'\in\Y^n : \L_r(\v y'|\v x)\leq L\},
\\ \nonumber
  \qmbox{where} L &:=& \L_r(\v y|\v x)
\eqa
For this to make sense, we have to assume (and will later assure)
that $\Rank_r(L)<\infty$, i.e.\ there are only
finitely many $\v y'\in\Y^n$ having loss smaller than $L$.
In a sense, $\rho=\Rank_r(\v y|\v
x)$ measures how compatible $\v y$ is with $r$; $\v y$ is the
$\rho$th most compatible with $r$.

Since the logarithm is a strictly monotone increasing function, we
can also consider the logarithmic rank $\LR_r(\v y|\v
x):=\log\Rank_r(\v y|\v x)$, which will be more convenient.

\begin{principle}[loss rank principle (LoRP) for classification]\label{prLRD}
For discrete $\Y$, the best classifier/regressor $r:\D\times\X\to\Y$
in some class $\R$ for data $D=(\v x,\v y)$ is the one of smallest
loss rank:
\beq\label{eqLRD}
  r^{best} \;=\; \arg\min_{r\in\R} \LR_r(\v y|\v x)
  \;\equiv\; \arg\min_{r\in\R} \Rank_r(\v y|\v x)
\eeq
where $\Rank_r$ is defined in \req{eqRank}.
\end{principle}

We give a simple example for which we can compute all ranks by hand
to help better grasping how the principle works, but the example is
too simplistic to allow any conclusion on whether the principle is
appropriate.

\fexample{exSD}{simple discrete}{
Consider $\X=\{1,2\}$, $\Y=\{0,1,2\}$, and two points
$D=\{(1,1),(2,2)\}$ lying on the diagonal $x=y$, with polynomial
(zero, constant, linear) least squares regression
$\R=\{r_0,r_1,r_2\}$ (see Ex.\ref{exPR}). $r_0$ is simply 0, $r_1$
the $y$-average, and $r_2$ the line through points $(1,y_1)$ and
$(2,y_2)$. This, together with the quadratic $\L$ for generic $\v
y'$ and observed $\v y=(1,2)$ (and fixed $\v x=(1,2)$), is
summarized in the following table
\beqn
\begin{array}{c|c|c|c}
  d & r_d(x|\v x,\v y') & \L_d(\v y'|\v x) & \L_d(D) \\ \hline
  0 &       0        & y'_1\!\,^2+y'_2\!\,^2 & 5         \\
  1 & \fr12(y'_1+y'_2)  & \fr12(y'_2-y'_1)^2  & \fr12       \\
  2 & (y'_2-y'_1)(x-1)+y'_1 & 0             & 0
\end{array}
\eeqn
From the $\L$ we can easily compute the Rank for all nine $\v
y'\in\{0,1,2\}^2$. Equal rank due to equal loss is
indicated by a $=$ in the table below. Whole equality groups are
actually assigned the rank of their right-most member, e.g.\ for
$d=1$ the ranks of $(y'_1,y'_2)=(0,1),(1,0),(2,1),(1,2)$ are all 7 (and not
4,5,6,7).
\beqn\def\trq{\hspace{4.2ex}}
\begin{array}{c|c|c}
    & \Rank_{r_d}(y'_1y'_2|12) \\
  d & \quad\trq 1 \trq 2 \trq 3 \trq 4 \trq 5 \trq 6 \trq 7 \trq 8 \trq 9 \!\! & \Rank_{r_d}(D)\\ \hline
  0 & y'_1 y'_2 = 00<01=10<11<02=20<21={\bf 12}<22                             & 8             \\
  1 & y'_1 y'_2 = 00=11=22<01=10=21={\bf 12}<02=20                             & 7             \\
  2 & y'_1 y'_2 = 00=01=02=10=11=20=21=22={\bf 12}                             & 9             \\
\end{array}
\eeqn
So LoRP selects $r_1$ as best regressor, since it has minimal rank
on $D$. $r_0$ fits $D$ too badly and $r_2$ is too flexible (perfectly
fits all $D'$).
} 

\paradot{LoRP for continuous $\Y$}
We now consider the case of continuous or measurable spaces $\Y$,
i.e.\ normal regression problems. We assume $\Y=\SetR$ in the
following exposition, but the idea and resulting principle hold for
more general measurable spaces like $\SetR^m$. We simply reduce the
model selection problem to the discrete case by considering the
discretized space $\Y_\eps=\eps\SetZ$ for small $\eps>0$ and
discretize $\v y\leadsto \v y_\eps\in\eps\SetZ^n$. Then
$\Rank_r^\eps(L):=\#\{\v y'_\eps\in\Y_\eps^n : \L_r(\v y'_\eps|\v
x)\leq L\}$ with $L=\L_r(\v y_\eps|\v x)$ counts the number of
$\eps$-grid points in the set
\beq\label{defV}
  V_r(L) \;:=\; \{\v y'\in\Y^n : \L_r(\v y'|\v x)\leq L\}
\eeq
which we assume (and later assure) to be finite, analogous to the
discrete case. Hence $\Rank_r^\eps(L)\cdot\eps^n$ is an
approximation of the {\em loss volume} $|V_r(L)|$ of set $V_r(L)$, and
typically
$\Rank_r^\eps(L)\cdot\eps^n=|V_r(L)|\cdot(1+O(\eps))\to|V_r(L)|$ for
$\eps\to 0$. Taking the logarithm we get $\LR_r^\eps(\v y|\v
x)=\log\Rank_r^\eps(L)=\log|V_r(L)|-n\log\eps +O(\eps)$. Since
$n\log\eps$ is independent of $r$, we can drop it in comparisons
like \req{eqLRD}. So for $\eps\to 0$ we can define the log-loss
``rank'' simply as the log-volume
\beq\label{defLRC}
  \LR_r(\v y|\v x) \;:=\; \log|V_r(L)|,
  \qmbox{where} L:=\L_r(\v y|\v x)
\eeq

\begin{principle}[loss rank principle for regression]\label{prLRC}
For measurable $\Y$, the best regressor $r:\D\times\X\to\Y$ in some
class $\R$ for data $D=(\v x,\v y)$ is the one of smallest loss
volume:
\beqn\label{eqLRC}
  r^{best} \;=\; \arg\min_{r\in\R} \LR_r(\v y|\v x)
  \;\equiv\; \arg\min_{r\in\R} |V_r(L)|
\eeqn
where $\LR$, $V_r$, and $L$ are defined in \req{defV} and \req{defLRC},
and $|V_r(L)|$ is the volume of $V_r(L)\subseteq\Y^n$.
\end{principle}

For discrete $\Y$ with counting measure we recover the discrete Loss Rank
Principle \ref{prLRD}.

\fexample{exSC}{simple continuous}{
Consider Example \ref{exSD} but with interval $\Y=[0,2]$.
The first table remains unchanged, while the second table becomes
\beqn
\begin{array}{c|c|c|c}
  d & V_d(L)=\{\v y'\in[0,2]^2: ...\} & |V_d(L)| & |V_d(\L_d(D))| \\ \hline
  0 & y'_1\!\,^2+y'_2\!\,^2\leq L  & {{2\sqrt{\max\{L-4,0\}\!}\;+}\atop {L({\pi\over 4}-\cos^{-1}(\min\{ {2\over\sqrt{L}},1\}))}} & \approx 3.6 \\
  1 & \fr12(y'_2-y'_1)^2\leq L      & 4\sqrt{2L}-2L & 3 \\
  2 & 0\leq L                      & 4 & 4
\end{array}
\eeqn
So LoRP again selects $r_1$ as best regressor, since it has smallest loss volume
on $D$.
} 

\paradot{Infinite rank or volume}
Often the loss rank/volume will be infinite, e.g.\ if we had chosen
$\Y=\SetZ$ in Ex.\ref{exSD} or $\Y=\SetR$ in Ex.\ref{exSC}. We will
encounter such infinities in Section \ref{secLM}. There are various
potential remedies. We could
modify (a) the regressor $r$ or %
(b) the $\L$ to make $\LR_r$ finite, %
(c) the Loss Rank Principle itself, or %
(d) find problem-specific solutions. %
Regressors $r$ with infinite rank might be rejected for
philosophical or pragmatic reasons. We will briefly consider (a) for
linear regression later, but to fiddle around with $r$ in a generic
(blackbox way) seems difficult. We have no good idea how to tinker
with LoRP (c), and also a patched LoRP may be less attractive. For
kNN on a grid we later use remedy (d). While in (decision) theory,
the application's goal determines the loss, in practice the loss is
often more determined by convenience or rules of thumb. So the $\L$
(b) seems the most inviting place to tinker with. A very simple
modification is to add a small penalty term to the loss.
\beq\label{eqLa}
   \L_r(\v y|\v x) \;\leadsto\;
   \L_r^\a(\v y|\v x) := \L_r(\v y|\v x)+\a||\v y||^2,
   \quad \a>0 \mbox{ ``small''}
\eeq
The Euclidian norm $||\v y||^2:=\sum_{i=1}^n y_i^2$ is default, but
other (non)norm regularizes are possible. The regularized
$\LR_r^\a(\v y|\v x)$ based on $\L_r^\a$ is always finite, since
$\{\v y:||\v y||^2\leq L\}$ has finite volume.
An alternative penalty $\a\v{\h y}^\trp\v{\h y}$, quadratic in
the regression estimates $\h y_i=r(x_i|\v x,\v y)$ is possible if
$r$ is unbounded in every $\v y\to\infty$ direction.

A scheme trying to determine a single (flexibility) parameter (like $d$
and $k$ in the above examples) would be of no use if it depended
on one (or more) other unknown parameters ($\a$), since varying through
the unknown parameter leads to any (non)desired result.
Since LoRP seeks the $r$ of smallest rank, it is natural to also
determine $\a$ by minimizing $\LR_r^\a$ w.r.t.\ $\a$. The good news
is that this leads to meaningful results.

\section{LoRP for Linear Models}\label{secLM}

In this section we consider the important class of linear regressors
with quadratic loss function. Since linearity is only assumed in $y$
and the dependence on $x$ can be arbitrary, this class is richer
than it may appear. It includes kNN (Example \ref{exKNN2}), kernel
(Example \ref{exKern}), and many other regressors. For linear
regression and $\Y=\SetR$, the loss rank is the volume of an
$n$-dimensional ellipsoid, which can efficiently be computed in time
$O(n^3)$ (Theorem \ref{thmLRL}). For the special case of projective
regression, e.g.\ linear basis function regression (Example
\ref{exLBFR}), we can even determine the regularization parameter
$\a$ analytically (Theorem \ref{thmLRP}).

\paradot{Linear regression}
We assume $\Y=\SetR$ in this section; generalization to $\SetR^m$ is
straightforward. A linear regressor $r$ can be written in the form
\beq\label{eqmj}
  \h y \;=\; r(x|\v x,\v y) \;=\; \sum_{j=1}^n m_j(x,\v x)y_j
   \quad\forall x\in\X
   \qmbox{and some} m_j:\X\times\X^n\to\SetR
\eeq
Particularly interesting is $r$ for $x=x_1,...,x_n$.
\beq\label{eqM}
  \h y_i \;=\; r(x_i|\v x,\v y) \;=\; \sum_j M_{ij}(\v x)y_j
  \qmbox{with} M:\X^n\to\SetR^{n\times n}
\eeq
where matrix $M_{ij}(\v x)=m_j(x_i,\v x)$. Since LoRP needs $r$ only on the
training data $\v x$, we only need $M$.

\fexample{exKNN2}{kNN ctd.}{
For kNN of Ex.\ref{exKNN} we have
$m_j(x,\v x)={1\over k}$ if $j\in\N_k(x)$ and 0 else, and
$M_{ij}(\v x)={1\over k}$ if $j\in\N_k(x_i)$ and 0 else.
} 

\fexample{exKern}{kernel regression}{
Kernel regression takes a weighted average over $\v y$,
where the weight of $y_j$ to $y$ is proportional to
the similarity of $x_j$ to $x$, measured by a kernel
$K(x,x_j)$, i.e.\ $m_j(x,\v x)=K(x,x_j)/\sum_{j=1}^n K(x,x_j)$.
For example the Gaussian kernel for $\X=\SetR^m$ is
$K(x,x_j)=\e^{-||x-x_j||_2^2/2\s^2}$.
} 

\fexample{exLBFR}{linear basis function regression, LBFR}{
Let $\phi_1(x),...,\phi_d(x)$ be a set or vector of ``basis''
functions often called ``features''. We place no restrictions on
$\X$ or $\v\phi:\X\to\SetR^d$. Consider the class of functions
linear in $\v\phi$:
\beqn
  \F_d \;=\; \{ f_{\v w}(x)=\textstyle{\sum_{a=1}^d} w_a\phi_a(x)=\v w^\trp\v\phi(x) : \v w\in\SetR^d \}
\eeqn
For instance, for $\X=\SetR$ and $\phi_a(x)=x^{a-1}$ we would recover
the polynomial regression Example \ref{exPR}.
For quadratic loss function $\L(y_i,\h y_i)=(y_i-\h y_i)^2$ we have
\beqn
  \L_{\v w}(\v y|\v\phi)
  \;:=\; \sum_{i=1}^n(y_i-f_{\v w}(x_i))^2
  \;=\; \v y^\trp\v y - 2\v y^\trp\v\Phi\v w + \v w^\trp B \v w
\eeqn
where matrix $\v\Phi$ is defined by $\v\Phi_{ia}=\phi_a(x_i)$ and
$B$ is a symmetric matrix with
$B_{ab}=\sum_{i=1}^n\phi_a(x_i)\phi_b(x_i)=[\v\Phi^\trp\v\Phi]_{ab}$.
The loss is quadratic in $\v w$ with minimum at $\v
w=B^{-1}\v\Phi^\trp\v y$. So the least squares regressor is $\h y=\v
y^\trp\v\Phi B^{-1}\v\phi(x)$, hence $m_j(x,\v x)=(\v\Phi
B^{-1}\v\phi(x))_j$ and $M(\v x)=\v\Phi B^{-1}\v\Phi^\trp$.
} 

Consider now a general linear regressor $M$ with quadratic loss
and quadratic penalty
\bqa\nonumber
  \L_M^\a(\v y|\v x) &=& \sum_{i=1}^n
  \left(y_i-\textstyle{\sum_{j=1}^n} M_{ij}y_j\right)^2+\a||\v y||^2
  \;=\; \v y^\trp S_\a\v y,
\\ \label{defSa}
  \qmbox{where\footnotemark} S_\a &=& (1\!\!1-M)^\trp(1\!\!1-M)+\a 1\!\!1
\eqa%
\footnotetext{The mentioned alternative penalty $\a||\v{\h
y}||^2$ would lead to
$S_\a = (1\!\!1-M)^\trp(1\!\!1-M)+\a M^\trp M$.
For LBFR, penalty $\a||\v{\h w}||^2$ is popular (ridge regression).
Apart from being limited to parametric regression, it
has the disadvantage of not being reparametrization invariant.
For instance, scaling $\phi_a(x)\leadsto\g_a\phi_a(x)$ doesn't
change the class $\F_d$, but changes the ridge regressor.}
($1\!\!1$ is the identity matrix). $S_\a$ is a symmetric matrix. For
$\a>0$ it is positive definite and for $\a=0$ positive semidefinite.
If $\l_1,...\l_n\geq 0$ are the eigenvalues of $S_0$, then
$\l_i+\a$ are the eigenvalues of $S_\a$. $V(L)=\{\v y'\in\SetR^n : \v
y'\!\,^\trp S_\a\v y'\leq L\}$ is an ellipsoid with the eigenvectors of
$S_\a$ being the main axes and $\sqrt{L/(\l_i+\a)}$ being their length.
Hence the volume is
\beqn
  |V(L)| \;=\; v_n\prod_{i=1}^n \sqrt{L\over \l_i+\a}
  \;=\; {v_n L^{n/2}\over\sqrt{\det S_\a}}
\eeqn
where $v_n=\pi^{n/2}/{n\over 2}!$ is the volume of the
$n$-dimensional unit sphere, $z!:=\G(z+1)$, and $\det$ is the
determinant. Taking the logarithm we get
\beq\label{eqLRL1}
  \LR_M^\a(\v y|\v x)
  \;=\; \log|V(\L_M^\a(\v y|\v x))|
  \;=\; \fr n2\log(\v y^\trp S_\a\v y)-\fr12\log\det S_\a +\log v_n
\eeq
Consider now a {\em class} of linear regressors $\R=\{M\}$,
e.g.\ the kNN regressors $\{M_k:k\in\SetN\}$ or
the $d$-dimensional linear basis function regressors $\{M_d:d\in\SetN_0\}$.

\begin{theorem}[LoRP for linear regression]\label{thmLRL}
For $\Y=\SetR$, the best linear regressor $M:\X^n\to\SetR^{n\times n}$
in some class $\M$ for data
$D=(\v x,\v y)$ is
\beq\label{eqLRL}
  M^{best} \;=\; \mathop{\arg\min}_{M\in\M,\a\geq 0}
  \{ \fr n2\log(\v y^\trp S_\a\v y)-\fr12\log\det S_\a \}
  \;=\; \mathop{\arg\min}_{M\in\M\;\a\geq 0}
  \Big\{ {\v y^\trp S_\a\v y\over(\det S_\a)^{1/n}} \Big\}
\eeq
where $S_\a$ is defined in \req{defSa}.
\end{theorem}

Since $v_n$ is independent of $\a$ and $M$ it was possible to drop
$v_n$. The last expression shows that linear LoRP minimizes the $\L$
times the geometric average of the squared axes lengths of ellipsoid
$V(1)$. Note that $M^{best}$ depends on $\v y$ unlike the $M\in\M$.

\paradot{Nullspace of $S_0$}
If $M$ has an eigenvalue 1, then $S_0=(1\!\!1-M)^\trp(1\!\!1-M)$ has a zero
eigenvalue and $\a>0$ is necessary, since $\det S_0=0$. Actually
this is true for most practical $M$. Nearly all linear regressors
are invariant under a constant shift of $\v y$, i.e.\
$r(y_i+c|D)=r(y_i|D)+c$, which implies that $M$ has eigenvector
$(1,...,1)^\trp$ with eigenvalue 1. This can easily be checked for kNN
(Ex.\ref{exKNN}), Kernel (Ex.\ref{exKern}), and LBFR
(Ex.\ref{exLBFR}). Such a generic 1-eigenvector effecting all $M\in\M$
could easily and maybe should be filtered out by considering only
the orthogonal space or dropping these $\l_i=0$ when computing $\det
S_0$. The 1-eigenvectors that depend on $M$ are the ones where we
really need a regularizer $\a>0$ for. For instance, $M_d$ in LBFR has $d$
eigenvalues 1, and $M_{\text{kNN}}$ has as many eigenvalues 1 as there are
disjoint components in the graph determined by the edges $M_{ij}>0$
In general we need to find the optimal $\a$ numerically.
If $M$ is a projection we can find $\a_{min}$ analytically.

\paradot{Projective regression}
Consider a projection matrix $M=P=P^2$ with $d=\tr P$ eigenvalues 1,
and $n-d$ zero eigenvalues. For instance, $M=\v\Phi B^{-1}\v\Phi^\trp$
of LBFR Ex.\ref{exLBFR} is such a matrix, since $M\v\Phi=\v\Phi$ and
$M\v\Psi=0$ for $\v\Psi$ such that $\v\Phi^\trp\v\Psi=0$. This implies
that $S_\a$ has $d$ eigenvalues $\a$ and $n-d$ eigenvalues $1+\a$.
Hence
\bqa\nonumber
  \det S_\a &=& \a^d(1+\a)^{n-d}, \qmbox{where}
  S_\a = S_0+\a 1\!\!1 = 1\!\!1-P+\a 1\!\!1
\\ \nonumber
  \v y^\trp S_\a\v y &=& (\r+\a)\v y^\trp\v y, \qmbox{where}
  \r:={\v y^\trp S_0\v y\over\v y^\trp\v y} = 1-{\v y^\trp P\v y\over\v y^\trp\v y}
\\ \label{eqLRPa}
  \Rightarrow\quad \LR_P^\a &=& \fr n2\log\v y^\trp\v y + \fr n2\log(\r+\a)-
  {\textstyle{d\over 2}}\log\a - {\textstyle{n-d\over 2}}\log(1+\a)
\eqa
The first term is independent of $\a$. Consider $1-\r>{d\over n}$,
the reasonable region in practice. Solving $\partial
\LR_P^\a/\partial\a = 0$ w.r.t.\ $\a$ we get a minimum at
$\a=\a_{min}:={\r d\over(1-\r)n-d}$. After some algebra we get
\beq\label{eqLRPamin}\textstyle
  \LR_P^{\a_{min}} \;=\; \fr n2\log\v y^\trp\v y - \fr n2 \KL({d\over n}||1-\r),
  \qmbox{where} \KL(p||q) \;=\; p\log{p\over q}+(1-p)\log{1-p\over 1-q}
\eeq
is the relative entropy or Kullback-Leibler divergence.
Minimizing $\LR_P^{\a_{min}}$ w.r.t.\ $M$ is equivalent to
maximizing $\KL({d\over n}||1-\r)$. This is an unusual task,
since one mostly encounters $D$ minimizations.
For fixed $d$, $\LR_P^{\a_{min}}$ is monotone increasing in $\r$.
Since $\L_P^\a\propto\r+\a$, LoRP suggests to minimize $\L$ for fixed
model dimension $d$.
For fixed $\r$, $\LR_P^{\a_{min}}$ is monotone increasing in $d$,
i.e.\ LoRP suggests to minimize model dimension $d$ for fixed $\L$.
Normally there is a tradeoff between minimizing $d$ and $\r$,
and LoRP suggests that the optimal choice is the one that maximizes $\KL$.

\begin{theorem}[LoRP for projective regression]\label{thmLRP}
The best projective regressor $P:\X^n\to\SetR^{n\times n}$
with $P=P^2$ in some projective class $\cal P$ for data
$D=(\v x,\v y)$ is
\beqn\label{eqLRP}
  P^{best} \;=\; \arg\max_{P\in\cal P} \;\textstyle \KL({\tr P(x)\over n}||{\v y^\trp P(x)\v y\over\v y^\trp\v y}),
  \qmbox{provided} {\tr P\over n} < {\v y^\trp P\v y\over\v y^\trp\v y}
\eeqn
\end{theorem}

\section{Comparison to Gaussian Bayesian Linear Regression}\label{secBayes}

We now consider linear basis function regression (LBFR) from a
Bayesian perspective with Gaussian noise and prior, and compare it
to LoRP. In addition to the noise model as in PML, one also has to
specify a prior. Bayesian model selection (BMS) proceeds by
selecting the model that has largest evidence. In the special case
of LBFR with Gaussian noise and prior and an ML-II estimate for the
noise variance, the expression for the evidence has a similar
structure as the expression of the loss rank.

\paradot{Gaussian Bayesian LBFR / MAP}
Recall from Sec.\ref{secLM} Ex.\ref{exLBFR} that $\F_d$ is the class
of functions $f_{\v w}(x)=\v w^\trp\v\phi(x)$ ($\v w\in\SetR^d$) that
are linear in feature vector $\v\phi$. Let
\beq\label{eqGauss}
  \Gauss_N(\v z|\v\mu,\Sigma) \;:=\;
  {\exp(-\fr12(\v z-\v\mu)^\trp\Sigma^{-1}(\v z-\v\mu))
   \over(2\pi)^{N/2}\sqrt{\det\Sigma}}
\eeq
denote a general $N$-dimensional Gaussian distribution with mean
$\v\mu$ and covariance matrix $\Sigma$.
We assume that observations $y$ are perturbed from $f_{\v w}(x)$ by
independent additive Gaussian noise with variance $\b^{-1}$ and zero
mean, i.e.\ the likelihood of $\v y$ under model $\v w$ is
$\P(\v y|\v w)=\Gauss_n(\v y|\v\Phi\v w,\b^{-1}1\!\!1)$,
where $\v\Phi_{ia}=\v\phi_a(x_i)$.
A Bayesian assumes a prior (before seeing $\v y$) distribution on
$\v w$. We assume a centered Gaussian with covariance matrix $(\a
C)^{-1}$, i.e.\
$\P(\v w)=\Gauss_d(\v w|\v 0,\a^{-1}C^{-1})$.
From the prior and the likelihood one can compute the evidence and the posterior
\bqa\label{eqEv}
  \mbox{Evidence:}\qquad\quad\; \P(\v y) &=& \int\P(\v y|\v w)\P(\v w)d\v w
                                      \;=\; \Gauss_n(\v y|\v 0,\b^{-1}S^{-1})
\\ \nonumber
  \mbox{Posterior:}\qquad \P(\v w|\v y) &=& \P(\v y|\v w)\P(\v w)/P(\v y)
                                       \;=\; \Gauss_d(\v w|\v{\h w},A^{-1})
\eqa
\beq\label{eqBAMS}
  B:=\v\Phi^\trp\v\Phi, \quad
  A:=\a C+\b B, \quad
  M:=\b\v\Phi A^{-1}\v\Phi^\trp, \quad
  S:=1\!\!1-M, \quad
\eeq
\beqn
  \v{\h w}:=\b A^{-1}\v\Phi^\trp\v y, \quad
  \v{\h y}:=\v\Phi\v{\h w}=M\v y
\eeqn
A standard Bayesian point estimate for $\v w$ for fixed $d$ is the
one that maximizes the posterior (MAP) (which in the Gaussian case
coincides with the mean)
$\v{\h w}=\arg\max_{\v w}\P(\v w|\v y)=\b A^{-1}\v\Phi^\trp\v y$.
For $\a\to 0$, MAP reduces to Maximum Likelihood (ML), which
in the Gaussian case coincides with the least squares regression of
Ex.\ref{exLBFR}. For $\a>0$, the regression matrix $M$ is not a
projection anymore.

\paradot{Bayesian model selection}
Consider now a family of models $\{\F_1,\F_2,...\}$. Here the $\F_d$
are the linear regressors with $d$ basis functions, but in general
they could be completely different model classes. All quantities in
the previous paragraph implicitly depend on the choice of $\F$,
which we now explicate with an index. In particular, the evidence
for model class $\F$ is $\P_\F(\v y)$.
Bayesian Model Selection (BMS) chooses the model class (here $d$) $\F$ of
highest evidence:
\beqn
  \F^{\text{BMS}} \;=\;\arg\max_\F\P_\F(\v y)
\eeqn
Once the model class $\F^{\text{BMS}}$ is determined, the MAP (or other)
regression function $f_{{\v w}_{\F^{\text{BMS}}}}$ or $M_{\F^{\text{BMS}}}$ are
chosen. The data variance $\b^{-1}$ may be known or estimated
from the data, $C$ is often chosen $1\!\!1$, and $\a$ has to be chosen
somehow. Note that while $\a\to 0$ leads to a reasonable
MAP=ML regressor for fixed $d$, this limit cannot be used for BMS.

\paradot{Comparison to LoRP}
Inserting \req{eqGauss} into \req{eqEv} and taking the
logarithm we see that BMS minimizes
\beq\label{eqLLG1}
  -\log\P_\F(\v y) \;=\; \fr\b2\v y^\trp S\v y - \fr12\log\det S - \fr n2\log\fr\b{2\pi}
\eeq
w.r.t.\ $\F$. Let us estimate $\b$ by ML: We assume a broad prior
$\a\ll\b$ so that $\b{\partial S\over\partial\b}=O({\a\over\b})$ can
be neglected. Then ${\partial\log\P_\F(\v y)\over\partial\b} =
\fr12\v y^\trp S\v y-{n\over 2\b}+O({\a\over\b}n)=0$ $\Leftrightarrow$
$\b\approx \h\b:=n/(\v y^\trp S\v y)$. Inserting $\h\b$
into \req{eqLLG1} we get
\beq\label{eqLLG2}
  -\log\P_\F(\v y) \;=\;
  \fr n2\log\v y^\trp S\v y - \fr12\log\det S - \fr n2\log\fr{n}{2\pi\e}
\eeq
Taking an improper prior $\P(\b)\propto\b^{-1}$ and integrating out
$\b$ leads for small $\a$ to a similar result. The last term in
\req{eqLLG2} is a constant independent of $\F$ and can be ignored.
The first two terms have the same structure as in linear LoRP
\req{eqLRL}, but the matrix $S$ is different.
In both cases, $\a$ act as regularizers, so we may minimize over
$\a$ in BMS like in LoRP. For $\a=0$ (which neither makes sense in
BMS nor in LoRP), $M$ in BMS coincides with $M$ of Ex.\ref{exLBFR},
but still the $S_0$ in LoRP is the square of the $S$ in BMS. For $\a>0$, $M$ of
BMS may be regarded as a regularized regressor as suggested in
Sec.\ref{secLoRP} (a), rather than a regularized loss function (b) used
in LoRP. Note also that BMS is limited to (semi)parametric regression,
i.e.\ does not cover the non-parametric kNN Ex.\ref{exKNN} and Kernel
Ex.\ref{exKern}, unlike LoRP.

Since $B$ only depends on $\v x$ (and not on $\v y$), and all $\P$
are implicitly conditioned on $\v x$, one could choose $C=B$. In
this case, $M=\g\v\Phi B^{-1}\v\Phi^\trp$, with $\g={\b\over\a+\b}<1$
for $\a>0$, is a simple multiplicative regularization of projection
$\v\Phi B^{-1}\v\Phi^\trp$, and \req{eqLLG2} coincides with
\req{eqLRPa} for suitable $\a$, apart from an irrelevant additive
constant, hence minimizing \req{eqLLG2} over $\a$
also leads to \req{eqLRPamin}.

\section{Comparison to other Model Selection Schemes}\label{secOMS}

In this section we give a brief introduction to Penalized Maximum
Likelihood (PML) for (semi)parametric regression, and its major
instantiations, the Akaike and the Bayesian Information Criterion
(AIC and BIC), and the Minimum Description Length (MDL) principle,
whose penalty terms are all proportional to the number of parameters
$d$. The {\em effective} number of parameters is often much smaller than
$d$, e.g.\ if there are soft constraints like in ridge regression. We
compare MacKay's \cite{MacKay:92} trace formula for Gaussian
Bayesian LBFR and Hastie's et al.\ \cite{Hastie:01} trace formula
for general linear regression with LoRP.

\paradot{Penalized ML (AIC, BIC, MDL)}
Consider a $d$-dimensional stochastic model class like the Gaussian
Bayesian linear regression example of Section \ref{secBayes}. Let
$\P_d(\v y|\v w)$ be the data likelihood under $d$-dimensional model
$\v w\in\SetR^d$. The maximum likelihood (ML) estimator for fixed
$d$ is
\beqn
  \v{\h w} \;=\; \arg\max_{\v w}\P_d(\v y|\v w)
  \;=\; \arg\min_{\v w}\{-\log \P_d(\v y|\v w)\}
\eeqn
Since $-\log\P_d(\v y|\v w)$ decreases with $d$, we
cannot find the model dimension by simply minimizing over $d$
(overfitting). Penalized ML adds a complexity term to get
reasonable results
\beqn
  \h d \;=\; \arg\min_d\{-\log \P_d(\v y|\v{\h w}) + \mbox{Penalty}(d) \}
\eeqn
The penalty introduces a tradeoff between the first and second
term with a minimum at $\h d<\infty$. Various penalties have been
suggested: The Akaike Information Criterion (AIC) \cite{Akaike:73}
uses $d$, the Bayesian Information Criterion (BIC) \cite{Schwarz:78}
and the (crude) Minimum Description Length (MDL) principle use $\fr
d2\log n$ \cite{Rissanen:78,Gruenwald:04} for Penalty$(d)$.
There are at least {\em three important conceptual differences} to LoRP:
\begin{itemize}\parskip=0ex\parsep=0ex\itemsep=0ex
\item In order to apply PML one needs to specify not only a class
of regression functions, but a full probabilistic model $\P_d(\v y|\v w)$,
\item PML ignores or at least does not tell how to incorporate
a potentially given loss-function,
\item PML (AIC,BIC,MDL) is mostly limited to (semi)parametric
models (with $d$ ``true'' parameters).
\end{itemize}

We discuss two approaches to the last item in the remainder of this
section: AIC, BIC, and MDL are not directly applicable %
(a) for non-parametric models like kNN or Kernel regression, or %
(b) if $d$ does not reflect the ``true'' complexity of the model.
For instance, ridge regression can work even for $d$ larger than
$n$, because a penalty pulls most parameters towards (but not
exactly to) zero. %
MacKay \cite{MacKay:92} suggests an expression for the effective
number of parameters $d_{e\!f\!f}$ as a substitute for $d$ in case
of (b), and Hastie et al.\ \cite{Hastie:01} more generally also for
(a).

\paradot{The trace penalty for parametric Gaussian LBFR}
We continue with the Gaussian Bayesian linear regression example
(see Section \ref{secBayes} for details and notation). Performing
the integration in \req{eqEv}, MacKay \cite[Eq.(21)]{MacKay:92}
derives the following expression for the Bayesian evidence for $C=1\!\!1$
\bqa\label{eqLLMK}
  -\log\P(\v y) &=&
  (\a \h E_W+\b \h E_D) + (\fr12\log\det A -\fr d2\log\a) - \fr n2\log\fr\b{2\pi}
\\ \nonumber
  \h E_D &=& \fr12||\v\Phi\v{\h w}-\v y||_2^2, \quad
  \h E_W =\fr12||\v{\h w}||_2^2
\eqa
(the first bracket in \req{eqLLMK} equals $\fr\b2\v y^\trp S\v y$ and
the second equals $-\fr12\log\det S$, cf.\ \req{eqLLG1}).
Minimizing \req{eqLLMK} w.r.t.\ $\a$ leads
to the following relation:
\beqn
  0 \;=\; \textstyle{-\partial\log\P(\v y)\over\partial\a} \;=\;
  \h E_W +\fr12\tr A^{-1}-\fr d{2\a} \qquad
 ({\partial\over\partial\a}\log\det A=\tr A^{-1})
\eeqn
He argues that
$\a||\v{\h w}||_2^2$ corresponds to the effective number of
parameters, hence
\beq\label{eqdeffMcK}
  d^{\text{McK}}_{e\!f\!f} \;:=\; \a||\v{\h w}||_2^2
  \;=\; 2\a \h E_W \;=\; d-\a\tr A^{-1}
\eeq

\paradot{The trace penalty for general linear models}
We now return to general linear regression $\v{\h y}=M(\v x)\v y$
\req{eqM}. LBFR is a special case of a projection matrix $M=M^2$
with rank $d=\tr M$ being the number of basis functions. $M$ leaves
$d$ directions untouched and projects all other $n-d$ directions to
zero. For general $M$, Hastie et al. \cite[Sec.5.4.1]{Hastie:01}
argue to regard a direction that is only somewhat shrunken, say by a
factor of $0<\b<1$, as a fractional parameter ($\b$ degrees of
freedom). If $\b_1,...,\b_n$ are the shrinkages = eigenvalues
of $M$, the effective number of parameters could be defined as
\cite[Sec.7.6]{Hastie:01}
\beqn
  d^{\text{HTF}}_{e\!f\!f} \;:=\; \sum_{i=1}^n \b_i \;=\; \tr M
\eeqn
which generalizes the relation $d=\tr M$ beyond projections.
For MacKay's $M$ \req{eqBAMS}, $\tr M=d-\tr A^{-1}$,
i.e.\ $d^{\text{HTF}}_{e\!f\!f}$ is consistent with and generalizes
$d^{\text{McK}}_{e\!f\!f}$.

\paradot{Problems}
Though nicely motivated, the trace formula is not without problems.
First, since for projections, $M=M^2$, one could equally well have
argued for $d^{\text{HTF}}_{e\!f\!f}=\tr M^2$. Second, for kNN we have $\tr
M=\fr nk$ (since $M$ is $\fr 1k$ on the diagonal), which does not
look unreasonable. Consider now kNN' where we average over the $k$
nearest neighbors {\em excluding} the closest neighbor. For
sufficiently smooth functions, kNN' for suitable $k$ is still a
reasonable regressor, but $\tr M=0$ (since $M$ is zero on the
diagonal). So $d^{\text{HTF}}_{e\!f\!f}=0$ for kNN', which makes no sense
and would lead one to always select the $k=1$ model.

\paradot{Relation to LoRP}
In the case of kNN', $\tr M^2$ would be a better estimate for the
effective dimension. In linear LoRP, $-\log\det S_\a$ serves as
complexity penalty. Ignoring the nullspace of
$S_0=(1\!\!1-M)^\trp(1\!\!1-M)$ \req{defSa}, we can Taylor expand
$-\fr12\log\det S_0$ in $M$
\beqn
  -\fr12\log\det S_0 \;=\; -\tr\log(1\!\!1\!-\!M)
  \;=\; \sum_{s=1}^\infty \fr1s\tr(M^s)
  \;=\; \tr M + \fr12\tr M^2 + ...
\eeqn
For BMS \req{eqLLG2} with $S=1\!\!1-M$ \req{eqBAMS} we get half of
this value. So the trace penalty may be regarded as a leading order
approximation to LoRP. The higher order terms prevent
peculiarities like in kNN'.

\section{Outlook}\label{secMisc}

So far we have only scratched at the surface of the Loss Rank
Principle. LoRP seems to be a promising principle with a lot of
potential, leading to a rich field. In the following we briefly
summarize miscellaneous considerations, which may be elaborated on
in the future: Experiments, Monte Carlo estimates for non-linear
LoRP, numerical approximation of $\det S_\a$, LoRP for
classification, self-consistent regression, explicit expressions for
kNN on a grid, loss function selection, and others.

\paradot{Experiments}
Preliminary experiments on selecting $k$ in kNN regression confirm
that LoRP selects a ``good'' $k$. (Even on artificial data we cannot
determine whether the ``right'' $k$ is selected, since kNN is not a
generative model). LoRP for LBFR seems to be consistent with rapid
convergence.

\paradot{Monte Carlo estimates for non-linear LoRP}
For non-linear regression we did not present an efficient algorithm
for the loss rank/volume $\LR_r(\v y|\v x)$. The high-dimensional
volume $|V_r(L)|$ \req{defV} may be computed by Monte Carlo
algorithms. Normally $V_r(L)$ constitutes a small part of $\Y^n$, and
uniform sampling over $\Y^n$ is not feasible. Instead one should
consider two competing regressors $r$ and $r'$ and compute $|V\cap
V'|/|V|$ and $|V\cap V'|/|V'|$ by uniformly sampling from $V$ and
$V'$ respectively e.g.\ with a Metropolis-type algorithm. Taking the
ratio we get $|V'|/|V|$ and hence the loss rank difference
$\LR_r-\LR_{r'}$, which is sufficient for LoRP. The usual tricks and
problems with sampling apply here too.

\paradot{Numerical approximation of $\det S_\a$}
Even for linear regression, a Monte Carlo algorithm may be faster
than the naive $O(n^3)$ algorithm \cite{Bai:96}. Often $M$ is a very
sparse matrix (like in kNN) or can be well approximated by a sparse
matrix (like for Kernel regression), which allows to approximate
$\det S_\a$, sometimes in linear time \cite{Reusken:02}.

\paradot{LoRP for classification}
A classification problem is or can be regarded as a regression problem
in which $\Y$ is finite. This implies that we need to compute (count)
$\LR_r$ for non-linear $r$ somehow, e.g.\ approximately by Monte Carlo.

\paradot{Self-consistent regression}
So far we have considered only ``on-data'' regression. LoRP only
depends on the regressor $r$ on data $D$ and not on
$x\not\in\{x_1,...,x_n\}$. One can construct canonical regressors
for off-data $x$ from regressors given only on-data in the following
way:
We add a virtual data point $(x,y)$ to $D$, where $x$ is the
off-data point of interest. If we knew $y$ we could estimate $\h
y=r(x|\{(x,y)\}\cup D)$, but we don't know $y$. But if we require
consistency, namely that $\h y=y$, we get a canonical
estimate for $\h y$.
First, this bootstrap may ease the specification of the regression
models, second, it is a canonical way for interpolation (LoRP
can't distinguish between $r$ that are identical on $D$), and third,
many standard regressors (kNN, Kernel, LBFR) are self-consistent in
the sense that they are canonical.

\paradot{Explicit expressions for kNN on a grid}
In order to get more insight into LoRP, a case that allows an
analytic solution is useful. For k nearest neighbors classification
with $\v x$ lying on a hypercube of the regular grid $\X=\SetZ^d$
one can derive explicit expressions for the loss rank as a function
of $k$, $n$, and $d$. For $n\gg k\gg 3^d$, the penalty
$-\fr12\log\det S$ is proportional to $\tr M$ with proportionality
constant decreasing from about 3.2 for $d=1$ to 1.5 for
$d\to\infty$.

\paradot{LoRP for hybrid model classes}
LoRP is not restricted to model classes indexed by a
single integral ``complexity'' parameter, but may be applied more
generally to selecting among some (typically discrete) class of
models/regressors. For instance, the class could contain kNN {\em
and} polynomial regressors, and LoRP selects the complexity {\em
and} type of regressor (non-parametric kNN versus parametric
polynomials).

\paradot{General additive loss}
Linear LoRP $\v{\h y}=M(\v x)\v y$ of Section \ref{secLM} can easily
be generalized from quadratic to $\rho$-norm $\L_M(\v y|\v
x)=||\v y-\v{\h y}||_\rho^p$ (any $p$). For $\a=0$, $\v y^\trp S_0\v y$ in
\req{eqLRL1} becomes $||\v y-\v{\h y}||_\rho^2$ and $v_\rho$ the
volume of the unit $d$-dimensional $\rho$-norm ``ball''. Useful
expressions for general additive $\L_N=\sum_i h(y_i-\h y_i)$
can also be derived. Regularization may be performed by
$M\leadsto\gamma M$ with optimization over $\gamma<1$.

\paradot{Loss-function selection}
In principle, the loss function should be part of the problem
specification, since it characterizes the ultimate goal.
In reality, though, having to specify the loss function can be a
nuisance.
We could interpret the regularized loss \req{eqLa} as a class of
loss functions parameterized by $\a$, and $\arg\min_\a\LR_r^a$ as a
loss function optimization or selection. This suggests to choose in
general the loss function that has minimal loss rank. This leads to
sensible results if the considered class of loss functions is not
too large (e.g.\ all $\rho$-norm losses in the previous paragraph).
So LoRP can be used not only for model selection, but also for
loss function selection.

\addcontentsline{toc}{section}{\refname}
\begin{small}

\end{small}

\end{document}